\documentclass[12pt,a4paper]{article}
\usepackage{amsmath,amssymb,amsthm} % 数学公式、符号、定理环境
\usepackage{geometry} % 页边距设置
\usepackage[numbers]{natbib} % 参考文献引用
\usepackage{graphicx} % 图片支持（预留）
\usepackage{hyperref} % 超链接（可选，期刊多支持）
\usepackage{tikz}

\theoremstyle{plain}
\newtheorem{theorem}{Theorem}[section]
\newtheorem{lemma}[theorem]{Lemma}

\newtheorem{prop}[theorem]{Proposition}
\theoremstyle{definition}

\newtheorem{claim}{Claim}
\newcommand{\exstar}{\text{ex}^*} % 彩虹图兰数符号
\newcommand{\floor}[1]{\left\lfloor #1 \right\rfloor} % 下取整
\newcommand{\ceil}[1]{\left\lceil #1 \right\rceil} % 上取整
\newcommand{\normG}[1]{\left\Vert #1 \right\Vert} % 图的边数符号
\newcommand{\Ft}{F_t} % t-扇图
\newcommand{\Wn}{W_n} % 轮图
\newcommand{\Cn}{C_n} % 圈图
\newcommand{\Kts}{K_{1,t,s}} % 完全三部图K_{1,t,s}
\newcommand{\Dn}{\Delta(n)} % Δ(n)函数
\newcommand{\ex}{\text{ex}}

\begin{document}

% 文章标题（英文，后续补充）
\title{Exact values of rainbow Turán numbers for fan graphs and even wheel graphs}

% 作者信息（预留，按需补充）
\author{Xinmin Hou, \quad Daoguang Xiong\\
\small $^{1}$School of Mathematical Sciences,\\
\small University of Science and Technology of China, Hefei 230026, Anhui, China\\
\small $^{2}$ Hefei National Laboratory,\\
\small University of Science and Technology of China, Hefei 230088, Anhui, China\\
\small Email: \texttt{xmhou@ustc.edu.cn (X. Hou)},  \texttt{pq377@mail.ustc.edu.cn (D. Xiong)}}

%\email{[pq377@mail.ustc.edu.cn]}

\maketitle
% 摘要（预留，后续补充英文翻译）
\begin{abstract}
%A proper edge coloring of a graph $G$ is an assignment of colors to all edges of $G$ such that no two adjacent edges share the same color. 
An edge-colored graph is called rainbow if all its edges have distinct colors. For a fixed graph $H$, the rainbow Turán number $\exstar(n, H)$ is  the maximum number of edges in a properly edge-colored graph with $n$ vertices that does not contain a rainbow subgraph isomorphic to $H$.
A $t$-fan $\Ft$ ($t \geq 2$) is a graph formed by $t$ triangles sharing a common vertex. A wheel graph $\Wn$ is constructed by connecting a new vertex $x$ to all vertices of a cycle $\Cn$ with $n$ vertices.  Keevash, Mubayi, Sudakov, Verstra\"ete~({\small Combin. Probab. Comput., 2007}) showed that 
\[
\exstar(n, F_t) \geq \floor{\frac{n^2}{4}} + (t-1)\floor{\frac{n}{2}} 
\]
by constructing extremal graphs for $n\not\equiv 2\pmod 4$.

In this paper, we propose such a method: for a graph $H$, the upper bound of $\exstar(n, H)$ can be analyzed using the rainbow Turán number of $H'$, where $H'$ is obtained by deleting one vertex from $H$. By applying this method, we determine the exact values of the rainbow Turán numbers for two families of graphs when the number of vertices $n$ is sufficiently large. Specifically, we obtain: (1)  when $n > 230t$,
\[
\exstar(n, \Ft) = \floor{\frac{n^2}{4}} + (t-1)\floor{\frac{n}{2}} - \Dn,
\]
where $\Dn = 1$ if $n \equiv 2 \pmod{4}$, and $\Dn = 0$ otherwise; (2) for graphs $H$ satisfying $W_{2t} \subset H \subset \Kts$ ($t \leq s$), when $n$ is sufficiently large,
\[
\exstar(n, H) = 
\begin{cases} 
\floor{\frac{n^2}{4}} + \floor{\frac{(t-1)n}{2}} - \Dn & \text{if } t \text{ is odd}, \\
\floor{\frac{n^2}{4}} + \floor{\frac{(t-1)n}{2}} & \text{if } t \text{ is even}.
\end{cases}
\]
\end{abstract}

% 关键词（预留，按需补充）

% -------------------------------------------------------------------
% Figures
% -------------------------------------------------------------------

% 正文框架（各章节预留内容插入位置）
%\maketitle

\section{Introduction}
The rainbow Turán problem, a generalization of the classical Turán problem, was first proposed by Keevash, Mubayi, Sudakov, Verstra\"ete~\cite{KMSV07}. The classical Turán problem, initiated by Tur\'an~\cite{T41}, focuses on finding the maximum number of edges in an $n$-vertex graph $G$ that does not contain a specific subgraph $H$, denoted as $\text{ex}(n, H)$. The celebrated Erd\H{o}s--Stone Theorem~\cite{ES46} determines $\text{ex}(n, H)$ asymptotically for every non-bipartite $H$. In contrast, the rainbow Turán problem considers the maximum number of edges in a properly edge-colored $n$-vertex graph $G$ that does not contain a rainbow subgraph isomorphic to $H$ (an edge-colored graph is called \emph{rainbow} if all its edges have distinct colors), denoted as $\exstar(n, H)$.

%Keevash et al.~\cite{KMSV07} proved that $ex(n,H)\leq\exstar(n,H)\leq ex(n,H)+o(n^2)$.Further,if H is \emph{colour-critical}(deleting some edge reduces the chromatic number) then $\exstar(n,H)=ex(n,H)$ for n sufficiently large.This indicates that the rainbow Turán number of a graph is closely related to its Turán number, yet there remain substantial differences between them, which motivates us to study the rainbow Turán numbers of graphs whose Turán numbers have already been determined.

{
As mentioned in \cite{KMSV07}, the rainbow Turán problem is motivated by two main directions. First, it connects to additive number theory: a $\mathcal{B}^*_k$-set in an abelian group can be used to construct a properly edge-colored bipartite graph with no rainbow $C_{2k}$, linking $\exstar(n, C_{2k})$ to the study of such sets. Second, it naturally generalizes both the Turán problem and the Ramsey problem, as illustrated by the Canonical Ramsey Theorem of Erd\H{o}s and Rado~\cite{ER50} and the subsequent works of Alon, Jiang, Miller, and Pritikin~\cite{AJMP03}.
}

Numerous studies have been conducted on the rainbow Turán problem, a selection of which is listed below.
\begin{itemize}
\item For matchings, LeSaulnier, Stocker, Wenger, West~\cite{LSWW10} showed  that an edge-colored graph with minimum color degree $\geq k$ contains a rainbow matching of size at least $\floor{\frac{k}{2}}$ (the \emph{color degree} of a vertex is the number of distinct colors of edges incident to it).  Johnston, Rombach~\cite{JR20} established that $\text{ex}^*(n, kK_2) = \ex(n, kK_2)= \binom{k-1}{2} + (k-1)(n-k+1)$.

 \item For the family of cycles $\mathcal{C}$, Keevash et al.~\cite{KMSV07} showed that $\exstar(n, \mathcal{C})=\Omega(n\log n)$. The upper bound was first given by Das, Lee, Sudakov~\cite{DLS13} and subsequently improved by Janzer~\cite{J23} and Tomon~\cite{T22}. The current best upper bound, due to  Kim, Lee, Liu, Tran~\cite{KLLT24}, is $\exstar(n, \mathcal{C})=O(n\log^2 n)$;

    \item For even cycle $C_{2k}$, Keevash et al.~\cite{KMSV07} showed that $\exstar(n, C_{2k})=\Omega(n^{1 + \frac{1}{k}})$ and conjectured that $\exstar(n, C_{2k}) = O(n^{1 + \frac{1}{k}})$. This conjecture was partial answered by Das, Lee, Sudakov~\cite{DLS13}, then  confirmed by Janzer~\cite{J23},  and later improved by Kim et al.~\cite{KLLT24}.

%    Verstraëte and Jacques \cite{JAN23} proved that $\exstar(n, C_{2k}) = O(n^{1 + \frac{1}{k}})$};
    
    \item For complete bipartite graphs $K_{t,s}$, Keevash et al.~\cite{KMSV07} showed that  $\exstar(n, K_{t,s}) = O(n^{2 - \frac{1}{t}})$.

    \item  For paths $P_{l}$ (paths with $l$ edges), Johnston, Rombach \cite{JR20}  showed that $\exstar(n,P_l)\geq \frac{nl}{2}+O(1)$ for $n\geq3$, Halfpap, Anastasia \cite{HA23} conjectured that this is the best lower bound, it has already been proved for n=3,4 by Johnston, Palmer, Sarkar \cite{JPS17} and $n=5$ by Halfpap, Anastasia \cite{HA23}. For the upper bound, Ergemlidze, Győri, and Methuku~\cite{EGM19} showed  that $\exstar(n, P_{l}) \leq \left( \frac{9l}{7} + 5 \right)n$.
\end{itemize}

Regarding the relationship between rainbow Turán numbers and classical Turán numbers, {Keevash et al.~\cite{KMSV07}} proved (using the Erd\H{o}s--Stone Theorem~\cite{ES46}) the following proposition. 
\begin{prop}[\cite{KMSV07}]
The rainbow Tur\'an number $\exstar(n, H)$ satisfies
 \[
\ex(n, H) \leq \exstar(n, H) \leq \text{ex}(n, H) + o(n^2).
\]
Furthermore, if $H$ is color-critical (there exists an edge in $H$ such that removing this edge reduces the chromatic number of $H$), then $\exstar(n, H) = \text{ex}(n, H)$ when $n$ is sufficiently large.   
\end{prop}
Keevash et al. also remarked that $\exstar(n, H)$ is not equal to $\ex(n, H)$ for a general non-bipartite graph $H$.  For example, for a $t$-fan $F_t$ (a graph that consists of $t$ triangles that all share exactly
one common vertex, called the center of $F_t$, see Figure \ref{Fig.1}),  Erd\H{o}s, F\"uredi, Gould, and Gunderson~\cite{EFGG95} proved  that for $n \geq 50t^2$,
\[
\text{ex}(n, F_t) = 
\begin{cases} 
\lfloor{\frac{n^2}{4}}\rfloor + t^2 - t, & \text{if } t \text{ is odd}, \\
\lfloor{\frac{n^2}{4}}\rfloor + (t-1)\frac{3t}{2}, & \text{if } t \text{ is even}.
\end{cases}
\]
Keevash, et al~\cite{KMSV07} showed that 
\[
\exstar(n, F_t) \geq \floor{\frac{n^2}{4}} + (t-1)\floor{\frac{n}{2}}(>\ex(n, F_t), \text{ when $n$ is large)} 
\]
by constructing extremal graphs for $n\not\equiv 2\pmod 4$.

\begin{figure}[htbp]
\centering

%--- Figure 1: Fan graph F_t (t independent triangles sharing only the center) ---
\begin{tikzpicture}[scale=1.2, vertex/.style={circle, fill=black, minimum size=5pt, inner sep=0pt}]\label{Fig.1}
  % center vertex
  \node[vertex, label=below:{$c$}] (c) at (0,-0.5) {};

  % --- First triangle (v_1, v_2) ---
  \node[vertex, label=left:{$v_1$}] (v1) at (-2.2,0.6) {};
  \node[vertex, label=above left:{$v_2$}] (v2) at (-1.4,1.6) {};
  \draw[thick] (c) -- (v1) (c) -- (v2) (v1) -- (v2);

  % --- Second triangle (v_3, v_4) ---
  \node[vertex, label=above:{$v_3$}] (v3) at (-0.6,1.9) {};
  \node[vertex, label=above:{$v_4$}] (v4) at (0.2,2.0) {};
  \draw[thick] (c) -- (v3) (c) -- (v4) (v3) -- (v4);

  % --- dots for generalization ---
  \node at (1.0,1.7) {\LARGE $\cdots$};

  % --- t-th triangle (v_{2t-1}, v_{2t}) ---
  \node[vertex, label=right:{$v_{2t}$}] (v2t_1) at (2.2,0.6) {};
  \node[vertex, label=above right:{$v_{2t-1}$}] (v2t) at (1.4,1.6) {};
  \draw[thick] (c) -- (v2t_1) (c) -- (v2t) (v2t_1) -- (v2t);

  % caption
  \node[anchor=north, yshift=-0.3cm] at (current bounding box.south) {{\bf Figure 1:} Fan graph $\Ft$ ($t$ triangles share only the common center $c$)};
\end{tikzpicture}
\end{figure}

In this paper, we improve this construction to establish the following lower bound:  % for the case $n \equiv 2 \pmod{4}$.
\[
\exstar(n, F_t) \geq \floor{\frac{n^2}{4}} + (t-1)\floor{\frac{n}{2}} - \Delta(n), 
\]
where
\[
\Delta(n) = 
\begin{cases} 
1 & \text{if } n \equiv 2 \pmod{4}, \\
0 & \text{otherwise}.
\end{cases}
\]

%Using Theorem \ref{thm:1.1}, we 
We further prove that this lower bound is also an upper bound when $n$ is large, thereby obtaining the exact value of the rainbow Turán number for $F_t$.
\begin{theorem}\label{thm:F_t}
For $t\ge 2$, $n > 230t$, 
\[
\exstar(n, F_t) = \floor{\frac{n^2}{4}} + (t-1)\floor{\frac{n}{2}} - \Delta(n).
\]

\end{theorem}

To prove Theorem~\ref{thm:F_t}, we use a heterochromatic-triangle-counting argument to establish the following key lemma, which may also be of independent interest.

\begin{lemma}[Key lemma]\label{lem:key}
Let $H'$ be the graph obtained by deleting one vertex $x$ from a given graph $H$. If the rainbow Turán number of $H'$ satisfies 
\[
\exstar(n, H') \leq M n^{1+s}
\]
for some constants $M$ and $0 \leq s \leq 1$, then the following hold:
\begin{itemize}
\item[(1).] Let $G$ be a properly  edge-colored  graph on  $n$ vertices. If  $G$ contains at least $\frac{8}{3}M |E(G)| n^s$ triangles, then $G$ contains a rainbow subgraph isomorphic to $H$.

\item[(2).] The rainbow Turán number of $H$ satisfies 
\[
\exstar(n, H) \leq \frac{1}{4}n^2 + 2M n^{1+s}.
\]
\end{itemize}
\end{lemma}

%In \cite{KMSV07}, to argue that $\exstar(n, H)$ is not always asymptotically equal to $\text{ex}(n, H)$ when $n$ is large, the authors considered the rainbow Turán number of the $t$-fan $F_t$ (a graph formed by $t$ triangles sharing a common vertex). Erd˝os P, F¨uredi Z, Gould R J, and Gunderson D S proved in \cite{2} that for $n \geq 50t$,
%\[
%\text{ex}(n, F_t) = 
%\begin{cases} 
%\floor{\frac{n^2}{4}} + t^2 - t & \text{if } t \text{ is odd}, \\
%\floor{\frac{n^2}{4}} + (t-1)\frac{3t}{2} & \text{if } t \text{ is even}.
%\end{cases}
%\]

%Let $f(n) = f(n, t) = \floor{\frac{n^2}{4}} + (t-1)\floor{\frac{n}{2}} - \Delta(n)$. Our key idea is to note that deleting the common vertex of the triangles in $F_t$ results in a $t$-matching; thus, the upper bound on the number of triangles in a graph $G$ without rainbow $F_t$ can be derived from Theorem \ref{thm:1.1}. If the number of edges in $G$ exceeds $f(n)$, we can use the upper bound on the number of triangles to show that a subgraph of $G$ has a structure approximating a balanced complete bipartite graph. Further analysis of this structure leads to a contradiction.

As a further application of Lemma~\ref{lem:key}, we determine the exact rainbow Turán number for certain graphs, including even wheels $W_{2t}$ and tripartite graphs $K_{1,s,t}$ ($t\leq s$). Here a {\it wheel} $W_{k}$ is obtained by taking a  $k$-cycle $C_{k}$ (i.e., a cycle with $k$ vertices) and adding a new vertex (the hub) connected to every vertex of the cycle (see Figure 2).
%, and  connecting a new center vertex to all vertices of .
\begin{figure}[htbp]
\centering
%\hfill
%--- Figure 2: Wheel graph W_{2t} (draw W_8 as example: hub + 8-cycle) ---
\begin{tikzpicture}[scale=1.2, vertex/.style={circle, fill=black, minimum size=5pt, inner sep=0pt}]
  % hub
  \node[vertex, label=below:{$x$}] (x) at (0,0) {};
  % 8 rim vertices in a circle
  \foreach \i [evaluate=\i as \angle using \i*45-90] in {1,...,8} {
    \node[vertex] (u\i) at (\angle:1.8) {};
  }
  % spokes: x to each rim vertex
  \foreach \i in {1,...,8} {
    \draw[thick] (x) -- (u\i);
  }
  % cycle edges
  \foreach \i in {1,...,7} {
    \draw[thick] (u\i) -- (u\the\numexpr\i+1\relax);
  }
 % \draw[thick] (u8) -- (u1);
  \draw[thick, dashed] (u8) -- (u1);
  % labels
  \node[below right] at (u1) {$u_1$};
  \node[above right] at (u2) {$u_2$};
  %\node[above right] at (u3) {$u_3$};
  \node[above left] at (u7) {$u_{k-1}$};
  \node[below left] at (u8) {$u_{k}$};
 % \node at (0.2,-1.6) {$\cdots$};
  % caption
  \node[anchor=north, yshift=-0.2cm] at (current bounding box.south) {{\bf Figure 2:} Wheel graph $W_k$ (hub $x$ connected to a $k$-cycle $C_k$)};
%\caption{Fig.2}
%\label{Fig.2}
\end{tikzpicture}

\end{figure}

\begin{theorem}\label{thm:1.3wheel}  
For graphs $H$ satisfying $W_{2t} \subseteq H \subseteq K_{1,s,t}$ ($t \leq s$), when $n$ is sufficiently large,
\[
\exstar(n, H) = 
\begin{cases} 
\lfloor{\frac{n^2}{4}}\rfloor + \lfloor{\frac{(t-1)n}{2}} \rfloor- \Delta(n), & \text{if } t \text{ is odd}, \\
\lfloor{\frac{n^2}{4}} \rfloor+ \lfloor{\frac{(t-1)n}{2}}\rfloor, & \text{if } t \text{ is even}.
\end{cases}
\]
\end{theorem}

%Let $g(n) = g(n, t)$ denote the right-hand side of the equation in Theorem \ref{thm:1.3}. Since $W_{2t} \subset H \subset K_{1,t,s}$, it follows that 
%\[
%\exstar(n, W_{2t}) \leq \exstar(n, H) \leq \exstar(n, K_{1,t,s}).
%\]
%To prove Theorem \ref{thm:1.3}, we only need to show that $\exstar(n, W_{2t}) \geq g(n)$ and $\exstar(n, K_{1,t,s}) \leq g(n)$. The construction for the former is similar to the lower bound proof of Theorem \ref{thm:1.2}, and the proof for the latter follows the same idea as the upper bound proof of Theorem \ref{thm:1.2}.

\subsection*{Notation}
Let $G$ be a simple graph with vertex set $V(G)$ and edge set $E(G)$.
we denote the number of vertices by $|G| = |V(G)|$ and the number of edges by $\normG{G} = |E(G)|$.
For a vertex  $x \in V(G)$, let $N_G(x)$ be the set of neighbors of $x$ in $G$, and $d_G(x) = |N_G(x)|$ is its degree. We may omit the subscript $G$ in the notation when the graph $G$ is clear from context.
 For $X, Y \subset V(G)$, let $G[X]$ denote the induced subgraph of $G$ on $X$. Define $e(X, Y)$ as the number of edges with one endpoint in $X$ and the other in $Y$. When $X=Y$, we have $e(X,Y)=\normG{G[X]}$, and for $v\in X$, write $d_X(v)$ for the degree of $v$ in $G[X]$.  For graphs $G$ and $H$, let $N(G, H)$ be the number of copies of $H$ in $G$.

Denote $$f(n, t) = \floor{\frac{n^2}{4}} + (t-1)\floor{\frac{n}{2}} - \Delta(n)$$ (abbreviated as $f(n)$ when $t$ is fixed),
and $$g(n, t)=\begin{cases} 
\lfloor{\frac{n^2}{4}}\rfloor + \lfloor{\frac{(t-1)n}{2}} \rfloor- \Delta(n), & \text{if } t \text{ is odd}, \\
\lfloor{\frac{n^2}{4}} \rfloor+ \lfloor{\frac{(t-1)n}{2}}\rfloor & \text{if } t \text{ is even}.
\end{cases}$$ (abbreviated as $g(n)$ when $t$ is fixed) throughout the paper.

The rest of this article is arranged as follows. In Section 2, we prove the key lemma. We give the lower bounds of $\exstar(n, F_t)$ and $\exstar(n, W_{2t})$ by constructing extremal graphs. The exact values of them are established in Sections 4 and 5. We give some remarks and discussions in the last section. 

\section{Proof of the key lemma}
In this section, we prove Lemma~\ref{lem:key} by first establishing two relevant lemmas.

\begin{lemma}\label{lem:2.1}
Consider $m$ triangles with a proper edge coloring (i.e., the three edges of each triangle have distinct colors). There exists a partition of the color set $C = C_1 \cup C_2$ (where $C_1 \cap C_2 = \emptyset$) such that the number of heterochromatic triangles is at least $\frac{3m}{4}$. A triangle is called \emph{heterochromatic} if it contains edges {with colors} from both $C_1$ and $C_2$; otherwise, it is called \emph{monochromatic}.
\end{lemma}
\begin{proof}
Consider the set $C$ of all colors appearing on the edges of these $m$ triangles. Partition $C$ randomly into $C_1 \cup C_2$ by independently assigning each color to $C_1$ or $C_2$ with probability $\frac{1}{2}$ each.

Fix any one of the $m$ triangles, say with edge colors $\{a, b, c\}$. Since the triangle is properly edge-colored, the three colors are distinct. The triangle is monochromatic under the random partition if and only if all three colors $\{a, b, c\}$ fall into the same part, i.e., either all in $C_1$ or all in $C_2$.
Therefore,
\[
\Pr[\text{triangle is monochromatic}] = \Pr[a,b,c \in C_1] + \Pr[a,b,c \in C_2] = \left(\frac{1}{2}\right)^3 + \left(\frac{1}{2}\right)^3 = \frac{1}{4}.
\]
Consequently,
\[
\Pr[\text{triangle is heterochromatic}] = 1 - \frac{1}{4} = \frac{3}{4}.
\]

Now let $X$ be the random variable counting the number of heterochromatic triangles under this random partition. For each triangle $T_i$, define the indicator $X_i = \mathbf{1}[T_i \text{ is heterochromatic}]$. Then $X = \sum_{i=1}^{m} X_i$, and by linearity of expectation:

\[
\mathbb{E}[X] = \sum_{i=1}^{m} \mathbb{E}[X_i] = \sum_{i=1}^{m} \frac{3}{4} = \frac{3m}{4}.
\]
Since the expected number of heterochromatic triangles is $\frac{3m}{4}$, there must exist at least one partition $C = C_1 \cup C_2$ for which $X \geq \frac{3m}{4}$. 
\end{proof}

\begin{lemma}\label{lem:2.2}
For a graph $G$ with $n$ vertices and $m$ edges, the number of triangles in $G$ is at least
\[
\frac{m(4m- n^2)}{3n}.
\]
\end{lemma}
\begin{proof}%[Proof of Lemma \ref{lem:2.2}]
For any edge $e = xy$, the number of triangles with $e$ as one edge is equal to the size of the intersection of the neighborhoods of $x$ and $y$, i.e.,
\[
|N(x) \cap N(y)| = |N(x)| + |N(y)| - |N(x) \cup N(y)|.
\]
Since $|N(x) \cup N(y)| \leq n$ (as $G$ has $n$ vertices), we have
\[
|N(x) \cap N(y)| \geq d(x) + d(y) - n.
\]
%where $d(x) = |N(x)|$ and $d(y) = |N(y)|$ denote the degrees of $x$ and $y$, respectively.
Summing over all edges $e \in E(G)$, we have 
\[
3 N(G, K_3) \geq \sum_{x \in G} d^2(x) - nm,
\]
where this inequality holds since for a vertex $x$, its degree $d(x)$ is counted exactly $d(x)$ times in the sum (once for each edge incident to $x$). 
By the Cauchy-Schwarz inequality, $\sum_{x \in G} d^2(x) \geq \frac{1}{n} \left( \sum_{x \in G} d(x) \right)^2$. Since $\sum_{x \in G} d(x) = 2m$ (the handshaking lemma), substituting this into the inequality gives
\[
3 N(G, K_3) \geq \frac{(2m)^2}{n} - n m = \frac{m(4m - n^2)}{n}.
\]
Consequently, we have 
\[
N(G, K_3)\geq \frac{m(4m- n^2)}{3n}.
\]
This completes the proof.
\end{proof}

\begin{proof}[Proof of Lemma~\ref{lem:key}]
Let $G$ be a properly edge-colored graph with $n$ vertices and $\normG{G} > \frac{1}{4}n^2 + 2M n^{1+s}$. We first show that $G$ contains a rainbow subgraph isomorphic to $H$.

By Lemma \ref{lem:2.2}, the number of triangles in $G$ satisfies
\[
N(G,K_3) > \normG{G} \cdot \frac{4\left( \frac{1}{4}n^2 + 2M n^{1+s} \right) - n^2}{3n} = \normG{G} \cdot \frac{8M n^{1+s}}{3n} = \frac{8}{3}M \normG{G} n^s.
\]
By Lemma \ref{lem:2.1}, there exists a partition of the color set $C$ into $C = C_1 \cup C_2$ ($C_1 \cap C_2 = \emptyset$) such that the number of heterochromatic triangles in $G$ is greater than:
\[
\frac{3}{4} \cdot \frac{8}{3}M \normG{G} n^s = 2M \normG{G} n^s.
\]
Define the following sets and notations:
\begin{itemize}
    \item $N_i(x)$: the set of vertices adjacent to $x$ via edges colored in $C_i$ for $i=1,2$;
    \item $\normG{G[N_i(x)]}_{3-i}$: the number of edges in the induced subgraph $G[N_i(x)]$ that are colored in $C_{3-i}$ for $i=1,2$.
\end{itemize}
Note that for each $i\in\{1,2\}$, $\sum_{x \in G} \normG{G[N_i(x)]}_{3-i}$ counts the number of heterochromatic triangles with exactly one $C_{3-i}$-colored edge  (each such triangle is counted once at the vertex opposite the  edge with color in $C_{3-i}$). 
%Similarly, $\sum_{x \in G} \normG{G[N_2(x)]}_1$ counts the number of heterochromatic triangles with exactly one edge colored in $C_1$. 
Thus
\[
\sum_{x \in G} \normG{G[N_1(x)]}_2 + \sum_{x \in G} \normG{G[N_2(x)]}_1 > 2M \normG{G} n^s= M n^s \left( \sum_{x \in G} d_1(x) + \sum_{x \in G} d_2(x) \right).
\]
This implies that there exists at least one vertex $x \in G$ such that either $\normG{G[N_1(x)]}_2 > M n^s d_1(x)$ or $\normG{G[N_2(x)]}_1 > M n^s d_2(x)$. Without loss of generality, assume $\normG{G[N_1(x)]}_2 > M n^s d_1(x)$.

Consider the subgraph $T$ of $G$ with vertex set $N_1(x)$ and edge set consisting of all $C_2$-colored edges in $G[N_1(x)]$. Then  $\normG{T} = \normG{G[N_1(x)]}_2 > M n^s d_1(x)>M |T|^{1+s}$, where the last inequality holds since $d_1(x) = |N_1(x)| = |T|<n$.
%, we have $\normG{T} > M |T|^{1+s}$. 
By the hypothesis of Lemma~\ref{lem:key}, $T$ contains a rainbow subgraph isomorphic to $H'$ (the graph obtained by deleting one vertex from $H$).

Note that all vertices in $T$ are adjacent to $x$ via $C_1$-colored edges. Thus, the rainbow $H'$ in $T$ and the vertex $x$ together form a subgraph isomorphic to $H$. We verify that this $H$ is rainbow:
first, edges colored with colors in $C_1$ (connecting $x$ to $T$) are distinct from edges colored with colors in $C_2$ (forming $H'$ in $T$);
second, all $C_1$-colored edges are incident to $x$, so they have distinct colors (due to proper edge coloring);
finally, all $C_2$-colored edges form a rainbow $H'$ and thus, they also have distinct colors.
%Thus, all edges of $H$ have distinct colors, meaning $H$ is rainbow. 
Therefore, the rainbow Turán number of $H$ satisfies
\[
\exstar(n, H) \leq \frac{1}{4}n^2 + 2M n^{1+s}.
\]
This completes the proof of Lemma~\ref{lem:key}.
\end{proof}

\section{Lower bounds of $\exstar(n, F_t)$  and $\exstar(n, W_{2t})$} % Theorems~\ref{thm:F_t} and \ref{thm:1.3wheel}}
%\vspace{5pt}
\noindent{\textbf{Construction of the Properly Edge-Colored Graphs $(G_F,\phi)$ and $(G_W, \phi)$}:}
Start with a properly edge-colored complete bipartite graph $K_{a,b}$ with parts $A$ ($|A| = a$) and $B$ ($|B| = b$).
%(the specific values of $a$ and $b$ will be determined later).
Choose \(t\ge 2\)  (the same \(t\) for both parts, and $a,b\geq2t-1$).
Let \(C_1, C_2, \dots, C_t\) be \(t\) completely new colors (not used in the original \(K_{a,b}\)).
For each part \(X \in \{A,B\}\),  do the following:
\begin{enumerate}
\item \textbf{Construction the edge-colored graph {\bf $(G_F,\phi)$}}
   \begin{itemize}
   \item If \(|X|\) is even, take \(t-1\) \emph{perfect} matchings.
Since a complete graph of even order admits a 1-factorization {(see, e.g., a classical result of Petersen~\cite{P91})}, it suffices to select t-1 of these $|X|-1(\geq t-1)$ edge-disjoint 1-factors.
   \item If \(|X|\) is odd, take \(t-1\) \emph{near-perfect} matchings (each leaves exactly one vertex unmatched) such that the unmatched vertices are all distinct.
%Let $x=2k+1$,and the vertex set $X = \{0, 1, 2, \dots, 2k\}$ with the cyclic group $\mathbb{Z}_{2k+1} = \mathbb{Z}/(2k+1)\mathbb{Z}$.For $0\leq i\leq t-2$,let $M_i=\{(i+j,i+j+k)\mid j=1,2,\dots ,k\}\pmod{2k+1}$ is a near-perfect matching with vertex i unmatched,this completes the construction.
Since a complete graph of odd order admits a decomposition into $(|X|-1)/2(\geq t-1)$ edge-disjoint Hamiltonian cycles, we first select t-1 of these cycles. We then choose one vertex from each selected Hamiltonian cycle such that all chosen vertices are pairwise distinct. Deleting the designated vertex from each Hamiltonian cycle yields a path with an odd number of edges, and it suffices to take alternating edges along each resulting path to obtain a near-perfect matching of X.
   \end{itemize}
   Assign color \(C_i\) to the edges of the \(i\)-th matching (\(i=1,\dots,t-1\)).

   \item \textbf{Construction the edge-colored graph $(G_W,\phi)$}

   Let \(U_X\) be the set of vertices in \(X\) that are not yet incident to \(t-1\) new edges.  By construction, \(|U_X|=t-1\) when \(|X|\) is odd, and \(U_X = \emptyset\) when \(|X|\) is even. If \(U_X \neq \emptyset\) and \(|U_X| \ge 2\), add a maximum matching on \(U_X\) and color it with color \(C_t\).
 \end{enumerate}
 \noindent{\bf Remark A:}   The construction of $(G_W,\phi)$ ensures that at most one vertex in $X$ has an {\it internal degree} (the number of new edges incident to the vertex in $X$) of $t-2$, while all other vertices in $X$ have an internal degree of $t-1$.
 Thus the total number of new edges added to $X$ is $\floor{\frac{(t-1)|X|}{2}}$. This construction ensures the resulting graph remains a properly edge-colored graph.

\begin{lemma}\label{LEM: lowerFt}
For $n\geq4t-1$, we have 
%{\color{red}(CONDITIONs?)}
\[
\exstar(n, F_t) \geq \floor{\frac{n^2}{4}} + (t-1)\floor{\frac{n}{2}} - \Dn = f(n, t),
\]
\end{lemma}
\begin{proof}
    We first claim that $(G_F,\phi)$ does not contain a rainbow $t$-fan $F_t$. Suppose, for contradiction, that $G_F$ contains a rainbow $F_t$. Since the complete bipartite graph $K_{a,b}$ has no triangles, the rainbow $F_t$ must contain at least $t$ edges from the newly added matchings (i.e., the non-bipartite edges). However, these newly added edges use only $(t-1)$ distinct colors. This contradicts the definition of a rainbow subgraph (all edges must have distinct colors). Thus, $G_F$ is a properly edge-colored graph without rainbow $F_t$.

Next, we calculate the number of edges in $G_F$. Clearly,
\[
\normG{G_F} = ab+ (t-1) \left( \floor{\frac{a}{2}}+ \floor{\frac{b}{2}} \right).
\]
We analyze two cases based on the value of $n \mod 4$.

\vspace{5pt}
\noindent{\bf Case 1:} $n \not\equiv 2 \pmod{4}$.\\
%When $n$ is not congruent to 2 modulo 4, 
Let  $|A| = \floor{\frac{n}{2}}$ or $|B| = \ceil{\frac{n}{2}}$, then $b\geq a\geq 2t-1$ and at most one of $a = \floor{\frac{n}{2}}$ or $b = \ceil{\frac{n}{2}}$ is odd,. For such $n$, we have
\[
\floor{\frac{1}{2}\floor{\frac{n}{2}}} + \floor{\frac{1}{2}\ceil{\frac{n}{2}}} = \floor{\frac{n}{2}}.
\]
Substituting this into the edge count of $G_F$ gives
\[
\normG{G_F} = \floor{\frac{n^2}{4}} + (t-1)\floor{\frac{n}{2}}.
\]
Therefore, %the lower bound of the rainbow Turán number satisfies:
\[
\exstar(n, F_t) \geq \floor{\frac{n^2}{4}} + (t-1)\floor{\frac{n}{2}}.
\]

\vspace{5pt}
\noindent{\bf Case 2:} $n \equiv 2 \pmod{4}$.\\
%When $n \equiv 2 \pmod{4}$, 
Let $n = 4k + 2$ ($k\geq t$) for some integer $k \geq 1$ and let  $a=2k$ and  $b=2k+2$, then $b\geq a\geq2t-1$.
%We consider the complete bipartite graph $K_{2k, 2k+2}$ (where $|A| = 2k$ and $|B| = 2k+2$) and construct the graph $G$ by adding matchings as described earlier. Calculating the number of edges in $G$, we have 
Then
\begin{eqnarray*}
\normG{G_F}& =& 2k(2k+2) + (t-1)k + (t-1)(k+1)\\
& =& (2k+1)^2 - 1 + (t-1)(2k+1)\\
&=&\floor{\frac{n^2}{4}} + (t-1)\floor{\frac{n}{2}}-1.
\end{eqnarray*}
% so the edge count simplifies to the above form. 
Therefore, 
\[
\exstar(n, F_t) \geq \floor{\frac{n^2}{4}} + (t-1)\floor{\frac{n}{2}} - 1.
\]

Combining these two cases, we obtain that
%and recalling that $\Dn = 1$ if $n \equiv 2 \pmod{4}$ and $\Dn = 0$ otherwise, we obtain the final lower bound for the rainbow Turán number of $F_t$:
\[
\exstar(n, F_t) \geq \floor{\frac{n^2}{4}} + (t-1)\floor{\frac{n}{2}} - \Dn = f(n, t).
\]
%where $f(n, t)$ is defined as in Section 1. 
%This completes the proof of the lower bound.
% 此处粘贴“3.1 下界的证明”部分的英文翻译内容
\end{proof}

\begin{lemma}\label{LEM: lowerWt}
For graphs $H$ satisfying $W_{2t} \subseteq H \subseteq K_{1,s,t}$ ($t \leq s$), $n\geq4t-1$, we have 
\[
\exstar(n, H) \geq g(n, t).
\]
\end{lemma}
\begin{proof}
%    \begin{claim}
We first claim that $G_W$ contains no rainbow $W_{2t}$.    
%\end{claim}
%\noindent\textbf{Verification of No Rainbow $W_{2t}$}:
%Let the resulting graph after adding edges be $G$. 
For clarity, we define the following:
\emph{cross edges} are the edges of the original complete bipartite graph, connecting vertices from different partitions $A$ and $B$;
\emph{internal edges} are the newly added edges, connecting vertices within the same partition $A$ or $B$.

%We claim that $G$ does not contain a rainbow $W_{2t}$. 
Suppose for contradiction that $G_W$ contains a rainbow $W_{2t}$. Let $x$ be the hub of this rainbow $W_{2t}$; removing $x$ from $W_{2t}$ yields a rainbow $2t$-cycle $C_{2t}$.

Without loss of generality, assume the hub $x \in A$. Let $r$ be the number of vertices of $C_{2t}$ that lie in partition $A$. Since all these $r$ vertices are adjacent to $x$ and $d_A(x)\leq t-1$, 
%the internal degree constraint (each vertex’s internal degree $\leq t-1$) 
we have $r < t$. Thus, there are at most $r$ internal edges between $x$ and  $C_{2t}$. Next, analyze the $2t$ edges of $C_{2t}$: each cross edge has one endpoint in $A$ (one of the $r$ vertices) and one in $B$. Since $C_{2t}$ is a cycle, each vertex in $A$ is adjacent to at most 2 edges of $C_{2t}$. Thus  the total number of cross edges in $C_{2t}$ is at most $2r$. Therefore, the number of internal edges in $C_{2t}$ is at least $2t - 2r$. In summary, the rainbow $W_{2t}$ contains at least $(2t - 2r) + r = 2t - r>t$ internal edges. However, the internal edges of $G_W$ are colored with only $t$ distinct colors. Therefore, $G_W$ does not contain a rainbow $W_{2t}$.

%\vspace{5pt}
%\noindent{\textbf{2. Calculation of the Number of Edges in $G$}}
Next, we calculate the number of edges in $G_W$.
To ensure that the graph $G$ has exactly $g(n)$ edges, we choose the sizes of parts $A$ and $B$ according to the residue of $n$ modulo 4.
If $n \equiv 2 \pmod{4}$ (i.e., $n = 4k + 2$ for some positive integer $k$),  then set $|A| = 2k + 2$ and $|B| = 2k$; 
Otherwise set $|A| = \floor{\frac{n}{2}}$ and $|B| = \ceil{\frac{n}{2}}$. We always have $b\geq a\geq 2t-1$.  Then
%(this ensures $|A| \cdot |B| = \floor{\frac{n^2}{4}}$, the maximum number of edges in a complete bipartite graph with $n$ vertices).
%The total number of edges in $G$ is the sum of cross edges and internal edges.
the number of cross edges in $G$ is $|A| \cdot |B| = \floor{\frac{n^2}{4}}-\Dn$.  Note that $a+b=n$ and at least one of $a$ and $b$ is even. 
One can check directly that the total number of internal edges added to $A$ and $B$ is 
$$\floor{\frac{(t-1)a}{2}} + \floor{\frac{(t-1)b}{2}}=\floor{\frac{(t-1)n}{2}}.$$
%$$\color{red}\floor{\frac{(t-1)a}{2}} + \floor{\frac{(t-1)b}{2}}=\begin{cases} 
%\floor{\frac{(t-1)n}{2}} - \Dn & \text{if } t \text{ is odd}, \\ 
%\floor{\frac{(t-1)n}{2}} & \text{if } t \text{ is even}. ???
%\end{cases}$$
%For $n \not\equiv 2 \pmod{4}$, the total number of cross edges is $\floor{\frac{(t-1)n}{2}}$; for $n \equiv 2 \pmod{4}$ and $t$ odd, this sum equals $\floor{\frac{(t-1)n}{2}} - 1$ (i.e., $\floor{\frac{(t-1)n}{2}} - \Dn$).
Therefore, the total number of edges in $G$ is exactly $g(n)$. Since $G$ is a properly edge-colored graph on $n$ vertices with $\normG{G} = g(n)$ that contains no rainbow $W_{2t}$. We conclude that for any $W_{2t} \subseteq H \subseteq K_{1,s,t}$
\[
\exstar(n, H)\ge \exstar(n, W_{2t}) \geq g(n).
\]

%For completeness, we clarify the definition of a bipartite partition used in the construction: A \textbf{bipartite partition} of a graph $G$ is a pair $(X, Y)$ such that $N(G) = X \cup Y$ and $X \cap Y = \emptyset$. For this partition:
% Edges with both endpoints in $X$ or both in $Y$ are called \textbf{internal edges};
%Edges with one endpoint in $X$ and one in $Y$ are called \textbf{cross edges};
%The \textbf{internal degree} of a vertex is the number of internal edges incident to it.
\end{proof}
\noindent{\bf Remark B:} The construction above shows that $g(n)$ is the maximum number of edges in an $n$-vertex graph $G$ that admits a bipartition $(X, Y)$ where every vertex has an internal degree of at most $t-1$. This property will be used in the subsequent upper bound proof.

\section{Proof of Theorem~\ref{thm:F_t}}
%\subsection{Proof of the Lower Bound}
%Consider an arbitrary properly edge-colored complete bipartite graph with two disjoint
%parts A and B, where $|A| = a$ and $|B| = b$ (the specific values of a and b will be determined
%later). We add $(t-1)$ maximum matchings of $A$ into $A$, and color all edges in each matching with a new color (i.e., a color not used in $K_{a,b}$). We perform the same operation on $B$, but use the same $(t-1)$ new colors (the ones added for $A$) to color the matchings added to $B$. This results in a properly edge-colored graph $(G,\phi)$.

%\subsection{Proof of the Upper Bound}
By Lemma~\ref{LEM: lowerFt}, we only need to prove that for fixed $t$ and $n > 230t$,
\[
\exstar(n, F_t) \leq \floor{\frac{n^2}{4}} + (t-1)\floor{\frac{n}{2}} - \Dn = f(n, t) = f(n).
\]
%where $f(n, t)$ is defined as in Section 1. 
Our proof uses contradiction: suppose $G$ is a properly edge-colored graph with $|G| = n > 230t$, $\normG{G} > f(n)$, and containing no rainbow $F_t$. 
{We derive a contradiction by establishing additional auxiliary lemmas and analyzing the structure of $G$.

Recall the result for rainbow matching established by LeSaulnier et al~\cite{LSWW10}.
\begin{lemma}[\cite{LSWW10}]\label{lem:matching}
An edge-colored graph with minimum color degree $\geq k$ contains a rainbow matching of size at least $\floor{\frac{k}{2}}$.
%(the \emph{color degree} of a vertex is the number of distinct colors of edges incident to it);
\end{lemma}

As an application of the above lemma, we give a properly edge-colored version by showing the following result. 
\begin{lemma}\label{lem:propermatching}
Let $G$ be a properly edge-colored graph with $|G| = n$. If $\normG{G} \geq 2nt$ for some positive integer $t$, then $G$ contains a rainbow matching of size $t$.
\end{lemma}
\begin{proof}%[Proof of Lemma \ref{lem:3.2}]
Since $\normG{G} \geq 2nt$, the average degree of $G$ is $d = \frac{2\normG{G}}{n} \geq 4t$. 
We may assume $\delta(G)\ge 2t$. Indeed, if there exists a vertex $v \in G$ with $d(v)< 2t$, then deleting $v$ and all its incident edges yields a graph $G'$ with average degree 
\[
d' > \frac{2(\normG{G} - 2t)}{n-1} \geq \frac{4nt - 4t}{n-1} = \frac{4t(n-1)}{n-1} = 4t.
\]
Repeating this deletion process eventually produces a subgraph of $G$ with minimum degree at least $2t$.

%LeSaulnier T D, Stocker C, Wenger P S et al. proved in \cite{2} that an edge-colored graph with minimum color degree $\geq k$ contains a rainbow matching of size at least $\floor{\frac{k}{2}}$. Here, "color degree" of a vertex refers to the number of distinct colors of edges incident to it. 
Note that for a properly edge-colored graph, the color degree of a vertex equals its degree.
%(since adjacent edges have distinct colors). 
By Lemma~\ref{lem:matching}, $G$ contains a rainbow matching of size at least $\floor{\frac{\delta(G)}{2}} \geq t$. This completes the proof.
\end{proof}

As an application of Lemma~\ref{lem:propermatching}, we continue have the following lemma . 
\begin{lemma}\label{lem:edge-triangle31}
Let $G$ be a properly edge-colored graph with $|G| = n$ and containing no rainbow $F_t$. Then

(1). $\normG{G} \leq \frac{1}{4}n^2 + 4nt$;

(2). The number of triangles in $G$ is at most $\frac{16t}{3}\normG{G}$.
\end{lemma}
\begin{proof}%[Proof of Theorem \ref{thm:3.1}]
%We first prove an auxiliary lemma:
%Returning to the proof of Theorem \ref{thm:3.1}, we 
%Apply Lemma~\ref{lem:key}. 
Let $H'$ be the graph obtained from  $F_t$ by deleting the common vertex of the triangles, i.e., $H'$ is a matching of size $t$. By Lemma \ref{lem:propermatching},  $\exstar(n, H') \leq 2nt$.
Apply Lemma~\ref{lem:key} to $H'$ with parameters $M = 2t$ and $s = 0$, and note that $G$ is $F_t$-free, we obtain that:\\ 
%in Theorem \ref{thm:1.1}. Substituting into Theorem \ref{thm:1.1}:
(1). The number of triangles in $G$ is at most $\frac{16t}{3}\normG{G}$;\\
(2). $\normG{G} \leq \frac{1}{4}n^2 + 2Mn^{1+s} = \frac{1}{4}n^2 + 4nt$.

%This completes the proof.
\end{proof}

The following lemma further characterizes the structure of an $n$-vertex $F_t$-free graph when $n$ is large.

\begin{lemma}\label{lem:subgraph3.3}
Let $G$ be a properly edge-colored graph with $|G| = n > 230t$, $\normG{G} > f(n)$, and containing no rainbow $F_t$. Then there exists a subgraph $G'$ of $G$ such that:
(a) $|G'| = m > 150t$;
(b) $\normG{G'} > f(m)$;
(c) The minimum degree $\delta(G') > 0.46m$.
\end{lemma}
\begin{proof}%[Proof of Theorem \ref{thm:3.3}]
We construct a sequence of subgraphs by deleting vertices one by one: $G = G_0 \supset G_1 \supset \dots \supset G_{n-m}$, where $G_{i+1}$ is obtained by deleting a vertex $v \in G_i$ if its degree in $G_i$ is $\leq 0.46|G_i|$.

We first show that the process stops when $m > 150t$. Since $G$ contains no rainbow $F_t$, $G_{n-m}$ also contains no rainbow $F_t$. By Lemma~\ref{lem:edge-triangle31}, $\normG{G_{n-m}} \leq \frac{1}{4}m^2 + 4mt$. Thus
\begin{equation}\label{eq:e(G)}
\frac{1}{4}m^2 + 4mt \geq \normG{G_{n-m}} = \normG{G} - \sum_{i=0}^{n-m-1} \left( \normG{G_i} - \normG{G_{i+1}} \right).
\end{equation}
Since we delete vertices with degree $\leq 0.46|G_i|$, the sum of deleted edges satisfies 
$$\sum_{i=0}^{n-m-1} \left( \normG{G_i} - \normG{G_{i+1}} \right) \leq \sum_{i=0}^{n-m-1} 0.46(n-i) = 0.23(n+m+1)(n-m).$$ Substituting $\normG{G} > f(n) \geq \frac{1}{4}n^2 - 1 + (t-1)\floor{\frac{n}{2}} - 1$ (since $\Dn \leq 1$) into inequality (\ref{eq:e(G)}), we have 
\[
\frac{1}{4}m^2 + 4mt \geq \left( \frac{1}{4}n^2 - 2 + (t-1)\floor{\frac{n}{2}} \right) - 0.23(n+m+1)(n-m)\geq  0.02n^2 + 0.23m^2 .
\]
Thus
\[0.02m^2 + 4mt > 0.02n^2> 0.02(230t)^2 = 1058t^2.
\]
Solving the quadratic inequality for $m$, we obtain
\[
 m > \frac{-4t + \sqrt{16t^2 + 4 \times 0.02 \times 1058t^2}}{2 \times 0.02} \approx 150.7t > 150t.
\]

Next, we show that $\normG{G_{n-m}} > f(m)$. For $k > 60$, the difference of $f(k)$ satisfies
\[
f(k) - f(k-1) \geq \left( \frac{k^2}{4} - \frac{(k-1)^2}{4} - 1 \right) - 1 = 0.5k - 2.25 > 0.46k.
\]
Since $\normG{G} > f(n)$ and each deletion removes at most $0.46|G_i|$ edges (from a vertex with degree $\leq 0.46|G_i|$), we can inductively prove $\normG{G_i} > f(|G_i|)$ for all $0 \leq i \leq n-m$. In particular, $\normG{G_{n-m}} > f(m)$.

Setting $G' = G_{n-m}$ completes the proof.
\end{proof}

We prove more structural results for graphs $G$ with rainbow $F_t\subseteq G$ or rainbow $F_t\not\subseteq G$.

\begin{lemma}\label{lem:neighborintersecting3.5.6}
Let $G$ be a properly edge-colored graph. Then the following hold:
\begin{itemize}
\item[(1)] If there exists a vertex $x \in V(G)$ such that the number of vertices $y \in N(x)$ with $|N(y) \cap N(x)| \geq 6t$ is at least $t$, then $G$ contains a rainbow $F_t$ with center $x$.

\item[(2)] If there exists a vertex $x \in V(G)$ and a subset $A \subset N(x)$ with $|A| > t$ such that $|N(y) \cap (N(x) \setminus A)| \geq 5t$ for all $y \in A$, then $G$ contains a rainbow $F_t$ with center $x$.
\end{itemize}
\end{lemma}

\begin{proof}%[Proof of Lemma \ref{lem:3.5}]
Take $t$ vertices $y_1, y_2, \dots, y_t \in N(x)$ (or $y_1, \dots, y_t \in A$ for case (2)) satisfying $|N(y_i) \cap N(x)| \geq 6t$ 
(or $|N(y_i) \cap (N(x) \setminus A)|\ge 5t$ for case (2)) for $1 \leq i \leq t$. 
If we can sequentially select vertices $z_1, z_2, \dots, z_t$ such that $z_i \in N(x) \cap N(y_i)$ (or $z_i \in N(y_i)\cap (N(x)\setminus A)$ for case (2)) for each $i$ and the triangles $xy_1z_1$, $xy_2z_2$, $\dots, xy_tz_t$ have pairwise distinct colors, then these $t$ triangles form a rainbow $F_t$  with center $x$. 
%Here, , all $z_i$ are distinct, and $z_i \neq y_j$ for all $i, j \in \{1, \dots, t\}$;  form .

We now verify such $z_i$'s exist. When selecting $z_{k+1} \in N(x) \cap N(y_{k+1})$ (or $z_{k+1} \in N(y_{k+1})\cap (N(x)\setminus A)$), the following constraints each eliminate at most a certain number of candidates:
\begin{itemize}
    \item[1.] the edge $xz_{k+1}$ must receive a color different from the colors of  $y_1z_1, \dots, y_kz_k$; this rules out at most $k$ vertices;
\item[2.] the edge $y_{k+1}z_{k+1}$ must receive a color different from the colors of edges $xy_1, \dots, xy_t$, $xz_1, \dots, xz_k$, and $y_1z_1, \dots, y_kz_k$, this rules out at most $t + 2k$ vertices;
\item[3.] $z_{k+1}$ must be distinct from $y_1, \dots, y_t$ and $z_1, \dots, z_k$, excluding  at most $t + k$ vertices  (in case (2),  $z_{k+1}$ are automatically distinct from $y_1, \dots, y_t$ since $z_{i}\in N(x)\setminus A$, and only $z_1, \dots, z_k$ need to be avoided. Consequently, at most $k$ vertices are excluded in this case).

\end{itemize}
Therefore, at most $k + (t + 2k) + (t + k) = 2t + 4k$ (or $k + (t + 2k) + k = t +4k$ for case (2)) vertices are removed from $|N(x) \cap N(y_{k+1})|$ (or $N(y_{k+1})\cap (N(x)\setminus A)$ for case (2)). 
Since $k \leq t-1$, we have $2t+4k\le 2t + 4(t-1) = 6t - 4 < 6t\le |N(x) \cap N(y_{k+1})|$ (or $t+4k\le t + 4(t-1) = 5t - 4 < 5t\le |N(y_{k+1})\cap (N(x)\setminus A)|$ for case (2)). 
Thus at least one admissible vertex remains, proving that a suitable $z_{k+1}$ exists. 

%(2). Take $t$ vertices $y_1, \dots, y_t \in A$ and select $z_i \in N(y_i) \cap (N(x) \setminus A)$ (for $1 \leq i \leq t$) using the same method as in Lemma \ref{lem:3.5}. The only difference is that $z_i \in N(x) \setminus A$, so constraint 3 (distinct from $y_1, \dots, y_t$) is no longer needed. The total number of excluded vertices for $z_{k+1}$ is $k + (t + 2k) + k = 2t + 3k < 5t \leq |N(y_{k+1}) \cap (N(x) \setminus A)|$ (since $k \leq t-1$ and $2t + 3(t-1) = 5t - 3 < 5t$). Thus, $z_{k+1}$ exists, and the triangles form a rainbow $F_t$.

This completes the proof. 
\end{proof}

\begin{lemma}\label{lem:edgecup3.7}
Let $G$ be a properly edge-colored graph with $\normG{G} > f(n)$ and containing no rainbow $F_t$. Then there exists an edge $e = xy \in E(G)$ such that $|N(x) \cup N(y)| \geq {n - 16t}$.
\end{lemma}
\begin{proof} %[Proof of Lemma \ref{lem:3.7}]
By Lemma~\ref{lem:edge-triangle31},  $N(G, K_3)\le \frac{16t}{3}\normG{G}$. For an edge $e = xy$, define $h(e) = n - |N(x) \cup N(y)|$ (the number of vertices not in $N(x) \cup N(y)$), and let $h$ be the average value of $h(e)$ over all $e \in E(G)$. The number of triangles containing the edge $e$  is $$|N(x) \cap N(y)| = d(x) + d(y) - |N(x) \cup N(y)| = d(x) + d(y) - (n - h(e)).$$
Summing over all edges $e \in E(G)$, we have 
\[
3\cdot N(G, K_3) \geq \sum_{x \in G} d^2(x) - n\normG{G} + \sum_{e \in E(G)} h(e).
\]
By the Cauchy-Schwarz inequality, $\sum_{x \in G} d^2(x) \geq \frac{1}{n}\left( \sum_{x \in G} d(x) \right)^2 = \frac{4\normG{G}^2}{n}$. Substituting $N(G, K_3)\leq \frac{16t}{3}\normG{G}$ into this inequality, we obtain
\[
3 \cdot \frac{16t}{3}\normG{G} \geq \frac{4\normG{G}^2}{n} - n\normG{G} + h\normG{G}.
\]
Dividing both sides by $\normG{G}$ (since $\normG{G} > 0$) yields that
\[
16t \geq \frac{4\normG{G}}{n} - n + h.
\]
Since $\normG{G} > f(n) \geq \frac{n^2}{4}$ (for $t \geq 2$ and $n > 0$), we have  $\frac{4\normG{G}}{n} - n > 0$. Consequently,  $h < 16t$. Therefore, there exists an edge $e = xy$ with $h(e) < 16t$, i.e., $|N(x) \cup N(y)| > n - 16t$. This completes the proof.
\end{proof}

We further obtain additional structural information about the $F_t$-free graph $G$ when its order is large.
\begin{lemma}\label{lem:3.4vertexpartition}
Let $G$ be a properly edge-colored graph with $|G| = n > 150t$, $\normG{G} > f(n)$, minimum degree $\delta(G) > 0.46n$, and containing no rainbow $F_t$. Then there exists a partition of $V(G)$ into $A_1 \cup A_2$ (i.e., $A_1 \cap A_2 = \emptyset$) such that:\\
(1). For each $i=1,2$ and every vertex $x \in A_i$, $|N(x) \cap A_{3-i}| > 0.3n$;\\
(2). $|A_i| > 0.46n - t$ for $i = 1, 2$.
\end{lemma}

\begin{proof}%[Proof of Theorem \ref{thm:3.4}]
%We first prove three auxiliary lemmas that will be used repeatedly:
%Returning to the proof of Theorem \ref{thm:3.4},
By Lemma \ref{lem:edgecup3.7}, we can choose an edge $e = x_0y_0$ with $|N(x_0) \cup N(y_0)| > n - 16t$. Hence there are at most $16t$ vertices outside of $N(x_0) \cup N(y_0)$ in $G$.

By Lemma~\ref{lem:neighborintersecting3.5.6} (1), $N(x_0)$ contains at most $t-1$ vertices $y$ with $|N(y) \cap N(x_0)| \geq 6t$; otherwise, $G$ has a rainbow $F_t$. Deleting these (at most $t-1$) vertices from $N(x_0)$ gives $A'_1 \subset N(x_0)$. Similarly, $N(y_0)$ also has at most $t-1$ vertices $z$ with $|N(z) \cap N(y_0)| \geq 6t$. Deleting these from $N(y_0)$ yields $A'_2 \subset N(y_0)$.

We claim $A'_1 \cap A'_2 = \emptyset$. Suppose for contradiction that there exists $z \in A'_1 \cap A'_2$. Then
\begin{eqnarray*}
 d(z) &=& |N(z) \cap N(x_0)| + |N(z) \cap N(y_0)| + |N(z) \cap (V(G) \setminus (N(x_0) \cup N(y_0)))| \\
 &<& 6t + 6t + 16t = 28t.   
\end{eqnarray*} Therefore, we have $0.46n < d(z)<28t$.
This contradicts to the assumption that $n>150t$. 
%But $n > 150t$, so $0.46n > 0.46 \times 150t = 69t > 28t$, a contradiction. Thus, $A'_1 \cap A'_2 = \emptyset$.

Moreover, $|A'_1| \geq |N(x_0)| - (t-1) > 0.46n - t$ and $|A'_2| \geq |N(y_0)| - (t-1) > 0.46n - t$. Hence the sizes of $A_1'$ and $A_2'$ satisfy the requirement of the statement (2).

Next, we show $A_1'$ and $A_2'$ satisfy the requirement of the statement (1). For any $z \in A'_1$,
\begin{eqnarray*}
|N(z) \cap A'_2|&=&|N(z)|-|N(z) \cap A'_1|-|N(z) \cap (V(G) \setminus (A'_1 \cup A'_2))|\\
                &>& 0.46n - 6t - 18t\\
                &=& 0.46n - 24t \\
                &>& 0.3n,  
\end{eqnarray*}
where the second inequality holds since $|N(z) \cap A'_1| \leq |N(z) \cap N(x)| < 6t$, and $|V(G) \setminus (A'_1 \cup A'_2)| \leq 16t + 2(t-1) < 18t$; the last inequality holds since $n > 150t$.
%Thus:
%\[
%|N(z) \cap A'_2| > 0.46n - 6t - 18t = 0.46n - 24t > 0.3n \quad (\text{since } 0.46 \times 150t - 24t = 45t > 0.3 \times 150t = %45t).
%\]
The same holds for $z \in A'_2$.
%, so $A'_1$ and $A'_2$ satisfy condition 1.

Finally, we extend $A'_1$ and $A'_2$ to be a partition of $V(G)$. The number of vertices outside of $A'_1 \cup A'_2$ is at most $16t + 2(t-1) <18t$. For each $z\notin A'_1 \cup A'_2$, we claim that there exists an $i \in \{1,2\}$ such that $|N(z) \cap A'_i| > 0.3n$, and thus we add $z$ to $A'_i$. Suppose for contradiction that $|N(z) \cap A'_1| \leq 0.3n$ and $|N(z) \cap A'_2| \leq 0.3n$. Set $A = N(z) \cap A'_1$. Then
\[
|A| > |N(z)| - |N(z) \cap A'_2| - |V(G)\setminus (A'_1 \cup A'_2)| > 0.46n - 0.3n - 18t > t.
\]
For any $u \in A$,
\begin{eqnarray*}
|N(u) \cap (N(z) \setminus A)| &\geq& |N(u) \cap (N(z) \setminus A'_1)| \\
                              &= &|(N(u) \setminus A'_1) \cap (N(z) \setminus A'_1)|\\
                              &\geq& |N(u) \setminus A'_1| + |N(z) \setminus A'_1| - |V(G)\setminus A'_1|\\
                             & >& (0.46n - 6t) + 0.16n - (0.54n + t)\\
                             &=& 0.08n - 7t \\
                             &>& 5t,
\end{eqnarray*}
where the forth inequality holds since $|N(u) \setminus A'_1| > 0.46n - 6t$, $|N(z) \setminus A'_1| > 0.46n - 0.3n = 0.16n$, and $|V(G) \setminus A'_1| < n - (0.46n - t) = 0.54n + t$.
%|N(u) \cap (N(z) \setminus A)| \geq |N(u) \cap (N(z) \setminus A'_1)| = |(N(u) \setminus A'_1) \cap (N(z) \setminus A'_1)|.
%\]
%By the principle of inclusion-exclusion:
%\[
%|(N(u) \setminus A'_1) \cap (N(z) \setminus A'_1)| \geq |N(u) \setminus A'_1| + |N(z) \setminus A'_1| - |G \setminus A'_1|.
%\]
%Since $|N(u) \setminus A'_1| > 0.46n - 6t$, $|N(z) \setminus A'_1| > 0.46n - 0.3n = 0.16n$, and $|G \setminus A'_1| < n - (0.46n %- t) = 0.54n + t$:
%\[
%|N(u) \cap (N(z) \setminus A)| > (0.46n - 6t) + 0.16n - (0.54n + t) = 0.08n - 7t > 5t \quad .
%\]
By Lemma \ref{lem:neighborintersecting3.5.6} (2), $G$ contains a rainbow $F_t$ with center $z$, a contradiction. 
Therefore, by suitably placing each vertex outside $A_1'\cup A_2'$ into either $A'_1$ or $A'_2$, we  obtain a desired partition $A_1, A_2$ of $V(G)$.
%satisfying both conditions of Theorem \ref{thm:3.4}. 
This completes the proof.
\end{proof}

\begin{lemma}\label{lem:3.8containF_t}
Let $G$ be a properly edge-colored graph with $|G| = n > 150t$, $\normG{G} > f(n)$, and minimum degree $\delta(G) > 0.46n$. Suppose there exists a partition $V(G) = A_1 \cup A_2$ such that:  \\ 
1. For each $i=1,2$ and every vertex $x \in A_i$, we have $|N(x) \cap A_{3-i}| > 0.3n$;\\
2. For each $i = 1, 2$, $|A_i| > 0.46n - t$.\\
Then $G$ contains a rainbow $F_t$.
\end{lemma}

\begin{proof}%[Proof of Theorem \ref{thm:3.8}]
For $i=1,2$ and a vertex $x \in A_i$, define the internal degree of $x$ in $A_i$ to be $d'(x) = |N(x) \cap A_i|$.
%(the "internal degree" of $x$ in $A_i$, for $i = 1, 2$). 
We consider two cases.

\vspace{5pt}
\noindent{\textbf{Case 1}: There exists a vertex $x \in V(G)$ with $d'(x) \geq t$.}

Without loss of generality, assume $x \in A_1$ and $y_1, y_2, \dots, y_t \in N(x) \cap A_1$. For each $y_i \in \{y_1, \dots, y_t\}$,
\begin{eqnarray*}
|N(x) \cap N(y_i)| &\geq& |N(x) \cap N(y_i) \cap A_2|\\
                 &\geq& |N(x) \cap A_2| + |N(y_i) \cap A_2| - |A_2|\\
                 &>& 0.3n + 0.3n - (0.54n + t)\\
                 &=& 0.06n - t\\
                 &>& 6t.
\end{eqnarray*}
%By the principle of inclusion-exclusion and conditions of Theorem \ref{thm:3.8}:
%\[
%|N(x) \cap N(y_i) \cap A_2| \geq |N(x) \cap A_2| + |N(y_i) \cap A_2| - |A_2|.
%\]
where the third inequality holds since $|N(x) \cap A_2| > 0.3n$, $|N(y_i) \cap A_2| > 0.3n$, and $|A_2| < n - (0.46n - t) = 0.54n + t$, and the last inequality holds since $n>150t$.
%\[
%|N(x) \cap N(y_i) \cap A_2| > 0.3n + 0.3n - (0.54n + t) = 0.06n - t > 6t \quad .
%\]
By Lemma \ref{lem:neighborintersecting3.5.6}, $G$ contains a rainbow $F_t$ with center $x$.

\vspace{5pt}
\noindent{\textbf{Case 2:} For all $x \in V(G)$, $d'(x) \leq t-1$.}

We first show there exists $i \in \{1,2\}$ such that $\normG{G[A_i]} > (t-1)\floor{\frac{|A_i|}{2}}$. Suppose for contradiction that $\normG{G[A_1]} \leq (t-1)\floor{\frac{|A_1|}{2}}$ and $\normG{G[A_2]} \leq (t-1)\floor{\frac{|A_2|}{2}}$. Then
\[
\normG{G} \leq |A_1||A_2| + \normG{G[A_1]} + \normG{G[A_2]} \leq \floor{\frac{n^2}{4}} + (t-1)\left( \floor{\frac{|A_1|}{2}} + \floor{\frac{|A_2|}{2}} \right).
\]
Since $\normG{G} > f(n) = \floor{\frac{n^2}{4}} + (t-1)\floor{\frac{n}{2}} - \Dn$, we have
\[
\floor{\frac{n^2}{4}} + (t-1)\left( \floor{\frac{|A_1|}{2}} + \floor{\frac{|A_2|}{2}} \right) > \floor{\frac{n^2}{4}} + (t-1)\floor{\frac{n}{2}} - \Dn.
\]
 The only way to avoid a contradiction is $\Dn = 1$ and $\floor{\frac{|A_1|}{2}} + \floor{\frac{|A_2|}{2}} = \floor{\frac{n}{2}}$.
However, $\Dn = 1$ implies $n \equiv 2 \pmod{4}$. Suppose $n = 4k + 2$. Then $\floor{\frac{n}{2}} = 2k + 1$. Given that $|A_1| + |A_2| = 4k + 2$, $|A_1||A_2| = \floor{\frac{n^2}{4}} = (2k+1)^2 - 1$ implies $|A_1| = |A_2| = 2k + 1$. But then
\[
\floor{\frac{|A_1|}{2}} + \floor{\frac{|A_2|}{2}} = k + k = 2k \neq 2k + 1 = \floor{\frac{n}{2}},
\]
a contradiction. Therefore, there exists an $i \in \{1,2\}$ with $\normG{G[A_i]} > (t-1)\floor{\frac{|A_i|}{2}}$; without loss of generality, assume $i = 2$.

Let $|A_1| = a$, $|A_2| = b$, and $N = |A_1||A_2| - e(A_1, A_2)$ (the number of non-edges between $A_1$ and $A_2$). By the handshaking lemma, the sum of internal degrees satisfies
\begin{eqnarray}\label{eqn:A1A2}
\frac{1}{2}(t-1)n\geq\frac{1}{2}\sum_{x \in V(G)} d'(x) &=& \normG{G[A_1]} + \normG{G[A_2]}\nonumber\\
                                                        &=& \normG{G} - e(A_1, A_2)\nonumber\\
                                                       &\geq& f(n)+1 - (|A_1||A_2| - N) \nonumber\\
                                                       &\geq&(t-1)\floor{\frac{n}{2}} + N \\
                                                       &\geq&(t-1) \frac{n-1}{2}+N\nonumber,   
\end{eqnarray}
where the forth inequality holds since $f(n) = \floor{\frac{n^2}{4}} + (t-1)\floor{\frac{n}{2}} - \Dn$ and $|A_1||A_2| \leq \floor{\frac{n^2}{4}}$.
%Since $\sum_{x \in G} d'(x) \leq n(t-1)$  and $\normG{G} > f(n)$, we have 
%\[
%\frac{1}{2}(t-1)n \geq \normG{G} - e(A_1, A_2) \geq f(n)+1 - (|A_1||A_2| - N).
%\]
%Substituting $f(n) = \floor{\frac{n^2}{4}} + (t-1)\floor{\frac{n}{2}} - \Dn$ and $|A_1||A_2| \leq \floor{\frac{n^2}{4}}$:
%\[
%\frac{1}{2}(t-1)n \geq(t-1)\floor{\frac{n}{2}} + N \geq(t-1) \frac{n-1}{2}+N.
%\]
Consequently, we have $N \leq \frac{t-1}{2}$.

Define  $C = \{x \in A_1 \mid e(x, A_2)<|A_2| \text{ and } d'(x) = t-1\}$ and $D = \{x \in A_1 \mid e(x, A_2)=|A_2| \text{ and } d'(x) = t-1\}$.
Let $c=|C|$. Then $c \leq N \leq \frac{t-1}{2}$. Let $d=|D|$. We claim that $d > 60t$. 
By the inequality (\ref{eqn:A1A2}), we have  $\normG{G[A_1]} + \normG{G[A_2]} \geq (t-1)\floor{\frac{n}{2}} $. Since  $\normG{G[A_2]} \leq (t-1)\frac{b}{2}$, we have
\[
\normG{G[A_1]} \ge (t-1)\left( \floor{\frac{n}{2}} - \frac{b}{2} \right) \geq (t-1)\frac{a-1}{2}.
\]
%The sum of the internal degrees in $A_1 \setminus C$ is:
Summing $d'(x)$ over all vertices of $A_1 \setminus C$, we obtain
\[
\sum_{x \in A_1 \setminus C} d'(x) = \sum_{x \in A_1} d'(x) - \sum_{x \in C} d'(x) \geq 2\normG{G[A_1]} - c(t-1) \ge (a-c-1)(t-1).
\]
On the other hand, $\sum_{x \in A_1 \setminus C} d'(x) = \sum_{x \in D} d'(x) + \sum_{x \in A_1 \setminus (C \cup D)} d'(x) \leq (t-1)d + (t-2)(a - c - d)$.  Consequently, we have %Equating the two bounds:
\[
(t-2)(a - c) + d \ge (a-c-1)(t-1)
\]
i.e., $d\ge a - c - (t-1)$. Since $a > 0.46n - t$, $c \leq \frac{t-1}{2}$, and $n > 150t$, we have 
\[
d \ge  0.46 \cdot 150t - t - \frac{t-1}{2} - (t-1) > 69t - 2.5t = 66.5t > 60t.
\]

We now show $G$ contains a rainbow $F_t$. 
%by analyzing two subcases based on edge colors in $A_2$.
Take a vertex \( x \in D \), let \( N(x) \cap A_1 = \{y_1, y_2, \dots, y_{t-1}\} \).
By the assumption that $\normG{G[A_2]} > (t-1)\floor{\frac{b}{2}}$ and each color in a properly edge-colored graph has at most $\floor{\frac{b}{2}}$ edges (no adjacent edges share the same color), $G[A_2]$ contains at least $t$ distinct colors. Thus there must exist a color \( c_0 \) that is different from the color of every edge \( xy_1, xy_2, \dots, xy_{t-1} \). Pick an edge \( uv \) in \( A_2 \) colored with \( c_0 \).

Remove all vertices \( z \) in \( A_2 \) satisfying the following conditions:
\begin{enumerate}
    \item \( z = u \) or \(  v \);
    \item the end $z$ of the edge \( xz \) with color \( c_0 \) (there is at most one such vertex);
    \item the ends $z$ of the edges $y_iz$ (for $1\le i\le t-1$) with colors different from those of the edges \( xu, xv, uv \)  (there are at most \( 3(t-1) \) such vertices $z$).
    %The color set of edges \( zy_1, zy_2, \dots, zy_{t-1} \) has a non-empty intersection with the color set of edges \( xu, xv, uv \) (there are at most \( 3(t-1) \) such vertices).
\end{enumerate}
Let $R$ be the set of  vertices remaining in $A_2$. Then %We then obtain a set \( E \subset A_1 \), with
\[
|R| \geq b - (2 + 1 + 3(t-1)) > 0.46n - t - 3t > 65t.
\]
Since $x\in D$, we have $N(x)\cap A_2=A_2$. Consequently, $N(x)\cap R=R$. Hence, for each \( y_i \) ($1\le i\le t-1$),
\[
|N(y_i)\cap N(x) \cap R| = |N(y_i) \cap R| \geq |R| - N > 65t - \frac{1}{2}(t-1) > 6(t-1).
\]
By Lemma \ref{lem:neighborintersecting3.5.6} (1), the subgraph $G[A_1\cup R]$ (also $G$) contains a rainbow $F_{t-1}$ centered at $x$, consisting of triangles \( \triangle xy_1z_1, \triangle xy_2z_2, \dots, \triangle xy_{t-1}z_{t-1} \) (with \( z_i \in R\)).
 By the construction of \( R \), all colors appearing in this fan are distinct from those  of the triangle \( \triangle xuv \). Therefore, adding \( \triangle xuv \) to this fan yields a rainbow \( F_t \).

This completes the proof.
\end{proof}

Now we are ready to prove Theorem~\ref{thm:F_t}.
\begin{proof}[Proof of Theorem~\ref{thm:F_t}:]
By contradiction, suppose $G$ is a properly edge-colored graph with $|G| = n > 230t$, $\normG{G} > f(n)$, and no rainbow $F_t$. By Lemma~\ref{lem:edge-triangle31}, $G$ has a bounded number of edges and triangles. Furthermore, by Lemma~\ref{lem:subgraph3.3}, $G$ contains a subgraph $G'$ such that $|G'| > 150t$, $\normG{G'} > f(|G'|)$, and $\delta(G') > 0.46|G'|$. Then applying  Lemma~\ref{lem:3.4vertexpartition} to $G'$, we obtain a partition $V(G')=A_1 \cup A_2$ satisfying conditions (1) and (2). Finally, by Lemma~\ref{lem:3.8containF_t}, $G'$ (and hence $G$) contains a rainbow $F_t$, a contradiction. Therefore, we have
\[
\exstar(n, F_t) \leq f(n) = \floor{\frac{n^2}{4}} + (t-1)\floor{\frac{n}{2}} - \Dn.
\]
Combining with the lower bound (Lemma~\ref{LEM: lowerFt}), we conclude that $\exstar(n, F_t) = f(n)$ for $n > 230t$. 

This completes the proof of Theorem \ref{thm:F_t}.
\end{proof}
}

\section{Proof of Theorem~\ref{thm:1.3wheel}}
Throughout this section, we assume $n$ is sufficiently large. By Lemma~\ref{LEM: lowerWt}, it suffices to prove $\exstar(n, K_{1,t,s}) \leq g(n)$ for $t \leq s$.
%, where $g(n)$ is defined in the lower bound proof (Section 4.1). 
%Following the proof idea of Theorem~\ref{thm:F_t}, we decompose the proof into three auxiliary theorems (Theorems 4.1–4.3). 
Let $ \alpha, \gamma, \delta$ be  sufficiently small positive constants with $0 < \alpha \ll \gamma \ll \delta$ (where "$\ll$" denotes "much smaller than").

Similar to the proof of Lemma~\ref{lem:subgraph3.3}, by greedily deleting vertices of small degrees, we can easily obtain the following lemma.
\begin{lemma}\label{lem:4.1subgraph}
For any positive integer $m$ and sufficiently small $\alpha > 0$, there exists a sufficiently large $n$ such that: If $G$ is a properly edge-colored graph on $n$ vertices with $\normG{G} > g(n)$ and containing  no rainbow $K_{1,t,s}$, then $G$ contains a subgraph $G'$ satisfying:\\
1. $|G'| = m$;\\
2. $\normG{G'} > g(m)$;\\
3. The minimum degree $\delta(G') > (0.5 - \alpha)m$.
\end{lemma}

{
In what follows, we will always assume such a minimum degree condition.
\begin{lemma}\label{lem:4.3disjointsubsets}
Let $\gamma > 0$ be sufficiently small. Suppose $\alpha, \beta > 0$ satisfy $\alpha < \frac{1}{6}$ and $\gamma > 12\sqrt{\frac{\beta}{1 - 6\alpha}}+\alpha $ (such $\alpha, \beta$ exist, e.g., by taking them sufficiently small). If $G$ is an $n$-vertex graph with $\delta(G) > (\frac 12 - \alpha)n$ and $N(G, K_3)<\beta n^3$, 
%fewer than $\beta n^3$ triangles 
then $G$ contains two disjoint subsets $B_1, B_2 \subset V(G)$ with the following property:\\
(1). $|B_i| > (\frac 12 - \gamma)n$ for $i = 1, 2$.\\
(2). For each $x \in B_i$ ($i = 1, 2$), $e(x, B_{3-i}) > (\frac 12 - 3\gamma)n$.
\end{lemma}
\begin{proof}%[Proof of Lemma \ref{lem:4.3}]
%The proof proceeds in four steps:
%\vspace{5pt}
%\noindent{\bf Step 1: Property of triangle edges.}\\
For edge $e\in E(G)$, let $t(e)$ be the number of triangles containing $e$ in $G$. 
\begin{claim}\label{cl:triangle}
Every triangle in $G$ has an edge $e$ such that $t(e)>\left( \frac{1}{6} - \alpha \right)n$.   
\end{claim}
For any triangle $\triangle xyz$, summing $t(xy) + t(xz) + t(yz)$ over all vertices $u \in G$ gives:
\begin{eqnarray*}
t(xy) + t(xz) + t(yz) &\geq&\sum_{u \in V(G)} \left( e(u, \{x,y,z\}) - 1 \right) \\
                      &= & d(x) + d(y) + d(z) - n\\
                      &>& 3\left(\frac 12 - \alpha\right)n - n \\
                      &=& \left(\frac 12 - 3\alpha\right)n.
\end{eqnarray*}
where the third inequality holds since $\delta(G)>(\frac 12 - \alpha)n$.
%$e(u, \{x,y,z\})$ is the number of edges between $u$ and $\{x,y,z\}$. Since $\delta(G)>(\frac 12 - \alpha)n$, we have
%\[
%d(x) + d(y) + d(z) - n > 3\left(\frac 12 - \alpha\right)n - n = \left(\frac 12 - 3\alpha\right)n.
%\]
Consequently,  at least one edge in $\triangle xyz$ satisfies $t(e) > \frac{1}{3}\cdot(\frac 12 - 3\alpha)n = \left( \frac{1}{6} - \alpha \right)n$.

%\paragraph{Step 2: Construct triangle-free subgraph $G_1$}
Delete all edges $e$ with $t(e) > \left( \frac{1}{6} - \alpha \right)n$ to obtain $G_1$. Since each such edge belongs to at least $\left( \frac{1}{6} - \alpha \right)n$ triangles, the number of deleted edges is at most $\frac{3\beta n^3}{(\frac{1}{6} - \alpha)n} = \frac{18\beta n^2}{1 - 6\alpha}$.
%(each such edge is in at least $\left( \frac{1}{6} - \alpha \right)n$ triangles). 
By Claim~\ref{cl:triangle}, $G_1$ contains no triangles.

%\paragraph{Step 3: Extract high-degree subgraph $G_2$}
Next, delete vertices $x \in V(G_1)$ with $d_{G_1}(x) < (\frac 12 - \varepsilon)n$ (where $\varepsilon = 6\sqrt{\frac{\beta}{1 - 6\alpha}}+\alpha$) to obtain $G_2$. For each such deleted vertex, at least $(\frac 12 - \alpha)n - (\frac 12- \varepsilon)n = (\varepsilon - \alpha)n$ edges were removed during the construction of $G_1$. Consequently, the number of such deleted vertices is less than  
\[
 \frac{2 \cdot \frac{18\beta n^2}{1 - 6\alpha}}{(\varepsilon - \alpha)n} = 6\sqrt{\frac{\beta }{1 - 6\alpha}}n.
\]
For each vertex $x \in V(G_2)$, $d_{G_2}(x) > (\frac12 - \varepsilon)n - 6\sqrt{\frac{\beta }{1 - 6\alpha}}n>{(\frac 12 - \gamma)n}$.

%\paragraph{Step 4: Build partition $B_1, B_2$}
Finally, we build $B_1$ and $B_2$.
Take an edge $e = xy \in E(G_2)$. Set $B_1 = N_{G_2}(x)$ and $B_2 = N_{G_2}(y)$. Then for each $i=1,2$, $|B_i| > (\frac 12-\gamma)n$ since $\delta(G_2)>(\frac 12-\gamma)n$.  Since $G_2$ is triangle-free, $B_1 \cap B_2 = \emptyset$. For each $z \in B_i$, $e(z, B_{3-i}) \geq d_G(z) - |G \setminus (B_1 \cup B_2)| > (\frac12 - \gamma)n - 2\gamma n = (\frac 12 - 3\gamma)n$.
\end{proof}

The following lemma asserts that if $s'$ is sufficiently large relative to $t$ and $s$, a proper edge-coloring of $K_{1,t,s'}$ contains a rainbow $K_{1,t,s}$.
\begin{lemma}\label{lem:rainbowk_1ts4.4}
Let $s' \geq t^2 + (s - 1)(t + 1)^2 + 1$. Then any proper edge-coloring of $K_{1,t,s'}$ contains a rainbow $K_{1,t,s}$.
\end{lemma}
\begin{proof}%[Proof of Lemma \ref{lem:4.4}]
Denote the three vertex parts of $K_{1,t,s'}$ by $\{x\}$ (called the center), $B$ (of size $t$), and $C$ (of size $s'$).
%Let $K_{1,t,s'}$ have three disjoint vertex sets: center $x$, $B$ (with $|B| = t$), and $C$ (with $|C| = s'$). 
Since $K_{1,t,s'}$ is properly edge-colored, the edges connecting $x$ and $B$ use $t$ distinct colors. 
Thus the number of vertices in $C$ adjacent to  $B$ via edges whose colors belong to the set of $t$ colors used by the edges from $x$ to $B$ is  at most $t^2$.
%(each vertex in $B$ has at most $t$ such neighbors). 
Delete these $t^2$ vertices from $C$ to obtain $C'$.

We now greedily select $z_1, z_2, \dots, z_s$ from $C'$ as follows. After choosing $z_i$, we delete from the remaining $C'$ all vertices that share any color on an edge from $z_i$ to $B \cup \{x\}$. This ensures that each subsequent vertex will have all its edge colors distinct from those of previously selected vertices. At each step, at most $(t + 1)^2$ vertices are removed, because each of the $t + 1$ vertices in $B \cup \{x\}$ is incident to at most $t + 1$ vertices in $C'$ via a given color. Since $|C'| \geq s' - t^2 \geq (s - 1)(t + 1)^2 + 1$, we can continue until select all $z_1, \dots, z_s$. Together with $x$ and $B$, these vertices form a rainbow $K_{1,t,s}$.
%uch neighbors) such that the colors of edges between $z_{k+1}$ and $B \cup \{x\}$ are distinct from the colors of edges between $\{z_1, \dots, z_k\}$ and $B$. For each $z_i$, we deleting vertices in $C'$ with overlapping colors removes at most $(t + 1)^2$ vertices (each of the $t + 1$ vertices in $B \cup \{x\}$ has at most $t + 1$ such neighbors). Since $|C'| \geq s' - t^2 \geq (s - 1)(t + 1)^2 + 1$, we can select all $z_1, \dots, z_s$. These vertices form a rainbow $K_{1,t,s}$ with $x$ and $B$.
\end{proof}

The following lemma asserts that if a properly edge-colored graph $G$ contains a large number of triangles, then $G$ contains a rainbow $K_{1,t,s}$.
\begin{lemma}\label{lem:rainbowK1ts4.5}
For any $\beta > 0$ and sufficiently large $n$, if an $n$-vertex properly edge-colored graph $G$ contains at least $\beta n^3$ triangles, then $G$ contains a rainbow $K_{1,t,s}$.
\end{lemma}
\begin{proof}%[Proof of Lemma \ref{lem:4.5}]
Similar as in the above lemma, let $\{x\}$ (called the center), $B$ (of size $t$), and $C$ (of size $s$) be the three vertex parts of $K_{1,t,s}$.
Deleting the center $x$ from $K_{1,t,s}$ yields  a complete bipartite graph $K_{t,s}$. Keevash et al.~\cite{KMSV07} proved that $\exstar(n, K_{t,s}) = O(n^{2 - \frac{1}{t}})$, so there exists a constant $M$ such that $\exstar(n, K_{t,s}) \leq M n^{2 - \frac{1}{t}}$.

For sufficiently large $n$, $\beta n^3 > \frac{8}{3}M n^{3 - \frac{1}{t}} > \frac{8}{3}M \normG{G} n^{1 - \frac{1}{t}}$. By Lemma~\ref{lem:key}
%Theorem \ref{thm:1.1} (Section 1), 
$G$ contains a rainbow $K_{1,t,s}$.
\end{proof}

Now we are ready to give a `stability lemma' for rainbow $K_{1,t,s}$-free graphs. 
\begin{lemma}\label{lem:4.2bipartition}
Let $G$ be a properly edge-colored graph on $n$ (sufficiently large) vertices with $\normG{G} > g(n)$ and  $\delta(G) > (\frac 12- \alpha)n$. Suppose $G$ contains no rainbow $K_{1,t,s}$ and constants $0 < \alpha \ll \gamma \ll \delta$. Then there exists a bipartition $(A_1, A_2)$ of $V(G)$ such that:\\
1. For $i = 1, 2$, we have $|A_i| > (\frac12 - \gamma)n$.\\
2. For each $i\in\{1,2\}$ and every vertex $x \in A_i$, we have $e(x, A_{3-i})> (\frac12 - \delta)n$.
\end{lemma}
\begin{proof}%[Proof of Theorem \ref{thm:4.2}]
%We first prove three auxiliary lemmas to support the partition construction:
%Returning to Theorem \ref{thm:4.2}, 
Since $0 < \alpha \ll \gamma$, we choose $\beta$ small enough to satisfy Lemma \ref{lem:4.3disjointsubsets}. By Lemma \ref{lem:rainbowK1ts4.5}, $G$ has fewer than $\beta n^3$ triangles (otherwise, $G$ contains a rainbow $K_{1,t,s}$, a contradiction). Thus, by Lemma \ref{lem:4.3disjointsubsets}, $G$ contains two disjoint subsets $B_1, B_2 \subset V(G)$ satisfying that  $|B_i| > (\frac12-\gamma)n$ and  $e(x, B_{3-i}) > (\frac12- 3\gamma)n$ for each $i\in\{1,2\}$ and every $x\in B_i$.

We expand $(B_1, B_2)$ to a bipartition $(A_1, A_2)$ of $V(G)$. For any $x \notin B_1 \cup B_2$, assume $e(x, B_1) \geq e(x, B_2)$ (the other case is symmetric). Then
\begin{eqnarray*}
e(x, B_1) &\geq& \frac{1}{2}\left( d(x) - |V(G)\setminus (B_1 \cup B_2)| \right)\\
          &>& \frac{1}{2}\left( \big(\frac12-\alpha\big)n - 2\gamma n \right)\\
          &=& \left( \frac{1}{4} - \frac{\alpha}{2} - \gamma \right)n > t.
\end{eqnarray*}
Take $t$ vertices $C \subset N(x) \cap B_1$. Since each vertex in $B_1$ has at most $4\gamma n$ non-neighbors in $B_2$, there are at most $4\gamma t n$ vertices in $B_2$ that are non-adjacent to $C$. 
If $e(x, B_2) - 4\gamma t n \ge s' = t^2 + (s - 1)(t + 1)^2 + 1$, then $G$ contains a copy of $K_{1,t,s'}$.
By Lemma \ref{lem:rainbowk_1ts4.4}, $G$ contains a rainbow $K_{1,t,s}$, a contradiction.
Therefore,
$e(x, B_2) - 4\gamma t n < s'$. Consequently,
\begin{eqnarray*}
e(x, B_1) &\geq& d(x) - |V(G)\setminus (B_1 \cup B_2)| - e(x, B_2)\\
          &>& \left(\frac12 - \alpha\right)n - 2\gamma n - s' - 4\gamma t n\\
          &>& \left(\frac12 - \delta\right)n.    
\end{eqnarray*}
Adding $x$ to $A_1$ preserves the partition conditions. Repeating this process for all such $x$ yields $A_1, A_2$ that satisfy the conditions of Lemma~\ref{lem:4.2bipartition}.
\end{proof}

We now yields a contradiction from the following lemma.

\begin{lemma}\label{lem:4.6contradiction}
Let $G$ be a properly edge-colored graph on $n$ (sufficiently large) vertices with  $\delta(G) > (\frac12 - \alpha)n$. Suppose $\normG{G} > g(n)$ and $V(G)$ has a bipartition $(A_1, A_2)$ such that:\\
1. For each $i\in\{1,2\}$ and every $x \in A_i$, we have $e(x, A_{3-i}) > (\frac12 - \delta)n$.\\
2. For each $i = 1, 2$, $|A_i| > (\frac 12 - \gamma)n$.\\
Then $G$ contains a rainbow $K_{1,t,s}$.
\end{lemma}
\begin{proof}%[Proof of Theorem \ref{thm:4.6}]
By Remark A,
%after the lower bound construction (Section 4.1), 
$g(n)$ is the maximum number of edges in an $n$-vertex graph with a bipartition where every vertex has internal degree $\leq t-1$. Since $\normG{G} > g(n)$, $G$ has at least one vertex $x$ with internal degree $\geq t$.
%(internal degree = number of edges within its partition).
Without loss of generality, assume $x \in A_1$. Suppose $y_1, \dots, y_t \in N(x) \cap A_1$.
%(so internal degree of $x$ is $\geq t$). 
Since $e(y, A_2) > (\frac12 - \delta)n$ for each vertex $y\in A_1$ and $|A_2| < (\frac12 + \gamma)n$, each vertex in $A_1$ has at most $(\gamma + \delta)n$ non-neighbors in $A_2$.
%(since $|N(x) \cap A_2| > (0.5 - \gamma)n$ and $|A_2| > (0.5 - \delta)n$), 
Thus $A_2$ has at most $(\gamma + \delta)(t + 1)n$ vertices non-adjacent to $\{x, y_1, \dots, y_t\}$.
Therefore, for sufficiently large $n$, the number of vertices in $A_2$ adjacent to all of $\{x, y_1, \dots, y_t\}$ is at least
\[
|A_2| - (\gamma + \delta)(t + 1)n > (\frac12 - \gamma)n - (\gamma + \delta)(t + 1)n \geq s' = t^2 + (s - 1)(t + 1)^2 + 1.
\]
Consequently, $G$ contains a properly edge-colored $K_{1,t,s'}$.
%with $x$ and $y_1, \dots, y_t$. 
By Lemma \ref{lem:rainbowk_1ts4.4}, this $K_{1,t,s'}$ (thus $G$) contains a rainbow $K_{1,t,s}$.
\end{proof}
}

Finally, for graphs $H$ satisfying $W_{2t} \subset H \subset K_{1,t,s}$, we have:\\ 
1. $\exstar(n, W_{2t}) \leq \exstar(n, H) \leq \exstar(n, K_{1,t,s})$\\
2. By Lemmas~\ref{lem:4.1subgraph}, \ref{lem:4.2bipartition} and \ref{lem:4.6contradiction}, $\exstar(n, K_{1,t,s}) \leq g(n)$.\\
3. By Lemma~\ref{LEM: lowerWt}, $\exstar(n, W_{2t}) \geq g(n)$.\\
Combining these gives $\exstar(n, H) = g(n)$ for sufficiently large $n$. This completes the proof of Theorem~\ref{thm:1.3wheel}.

\section{Discussions and remarks}
In this paper, we propose the heterochromatic-triangle-counting argument  for estimating the upper bound of rainbow Turán numbers, specifically Theorems~\ref{thm:F_t} and \ref{thm:1.3wheel}. This method enables us to analyze the upper bound of $\exstar(n, H)$ using the rainbow Turán number of $H'$ (the graph obtained by deleting one vertex from $H$). Additionally, we derive an upper bound on the number of triangles in a graph that does not contain a rainbow $H$. By leveraging these two upper bounds, we can further constrain the structure of graphs without rainbow $H$, thereby achieving more precise estimates of the upper bound of $\exstar(n, H)$ for specific graphs $H$.

%Following this approach, we investigate the rainbow Turán numbers for two families of graphs and obtain their exact values when the number of vertices $n$ is sufficiently large.

%1. For the $t$-fan graph $F_t$ (a graph formed by $t$ triangles sharing a common vertex): When $n > 230t$, the exact value of its rainbow Turán number is
%   \[
%   \exstar(n, F_t) = \floor{\frac{n^2}{4}} + (t-1)\floor{\frac{n}{2}} - \Dn,
%   \]
%   where $\Dn$ is defined as $\Dn = 1$ if $n \equiv 2 \pmod{4}$, and $\Dn = 0$ otherwise.

%2. For graphs $H$ satisfying $W_{2t} \subset H \subset K_{1,s,t}$ ($t \leq s$, where $W_{2t}$ denotes the even wheel graph and $K_{1,s,t}$ denotes the complete 3-partite graph): When $n$ is sufficiently large, the exact value of their rainbow Turán number is
%   \[
%   \exstar(n, H) = 
%   \begin{cases} 
%   \floor{\frac{n^2}{4}} + \floor{\frac{(t-1)n}{2}} - \Dn & \text{if } t \text{ is odd}, \\
%   \floor{\frac{n^2}{4}} + \floor{\frac{(t-1)n}{2}} & \text{if } t \text{ is even}.
%   \end{cases}
%   \]

Below are supplementary explanations regarding the estimation of rainbow Turán numbers using Lemma~\ref{lem:key}:
\begin{itemize}
    \item In Lemma~\ref{lem:key}, if $s < 1$ and $H$ is not  bipartite, then the rainbow Turán number of $H$ satisfies
   \[
   \exstar(n, H) = \frac{n^2}{4} + O(n^{1+s}).
   \]
This result is relatively satisfactory, as it shows that $\exstar(n, H)$ and $\ex(n, H)$ are asymptotically of the same order. Furthermore, using the constraint on the number of triangles and following the proof ideas of Theorems \ref{thm:F_t} and \ref{thm:1.3wheel}, it is easy to prove that a graph $G$ with $\normG{G} \geq \frac{n^2}{4}+o(n^2)$ and containing no rainbow $H$ must have a structure approximating a balanced complete bipartite graph. This structural insight can be used to obtain more precise estimates, making Lemma~\ref{lem:key} suitable for estimating rainbow Turán numbers of such graphs. For example, the suspension graph $\widehat{P_t}$ of a path $P_t$ is the graph obtained by join of a new vertex to every vertex of $P_t$. Using the proof structure developed in this paper, we can obtain $$\exstar(n,\widehat{P_t})=f(n,\lfloor \frac t2\rfloor),\, (t\geq 2).$$

% {\color{red}$\exstar(n,\widehat{P_{2t}})=\exstar(n,\widehat{P_{2t+1}})=f(n,t)\, (t\geq2)$. }

%\item Lemma~\ref{lem:key} allows us to iteratively estimate the upper bound of $\exstar(n, H)$ by reducing $H$ to graphs $H'=H-x_1, H''=H'-x_2, \dots$, but this process yields very crude bounds. For instance, the resulting upper  bound of  $\exstar(n, H)$ is at least $\frac{3n^2}4$ through one step iteration.
%into the estimation of the upper bound of $\exstar(n, H')$. Theoretically, we could further delete a vertex from $H'$ to obtain $H''$, and then use Lemma~\ref{lem:key} and the rainbow Turán number of $H''$ to estimate the upper bound of $\exstar(n, H')$. However, this iterative deletion approach yields very rough results: the upper bound of $\exstar(n, H')$ estimated in this way is at least $\frac{n^2}{4}$, and when this bound is used to estimate $\exstar(n, H)$, the resulting upper bound is at least $\frac{n^2}{4} + 2 \cdot \frac{n^2}{4} = \frac{3n^2}{4}$—a value very close to the number of edges in a complete graph $K_n$. Further deleting vertices from $H''$ to obtain $H'''$ and applying Lemma~\ref{lem:key} is infeasible, as the final estimated upper bound will exceed the number of edges in $K_n$, rendering the result meaningless.

\item When applying Lemma~\ref{lem:key}, the coefficient $2M$ in the conclusion
   \[
   \exstar(n, H) \leq \frac{1}{4}n^2 + 2M n^{1+s}
   \]
   can be improved by refining Lemma \ref{lem:2.1} (on heterochromatic triangles) and Lemma \ref{lem:2.2} (on the lower bound of the number of triangles). For example, the estimation of the number of triangles in Lemma \ref{lem:2.2} (i.e., a graph $G$ contains at least $\normG{G} \cdot \frac{4\normG{G} - n^2}{3n}$ triangles) can be replaced by more precise results under specific constraints, such as:\\
   - Bollobás~\cite{B76} proved that if an $n$-vertex graph $G$ satisfies $\frac{n^2}{4} \leq \normG{G} \leq \frac{n^2}{3}$, then $G$ contains at least $\frac{n}{9}(4\normG{G} - n^2)$ triangles;\\
   - Lovász and Simonovits~\cite{LS83} proved that if an $n$-vertex graph $G$ satisfies $\normG{G} = \floor{\frac{n^2}{4}} + t$ (where $t < \frac{n}{2}$), then $G$ contains at least $t\floor{\frac{n}{2}}$ triangles. \\
Depending on the specific properties of the graph $H$ or the constraints on the number of edges of $G$, appropriate theorems can be selected to estimate the number of triangles, thereby refining the upper bound of $\exstar(n, H)$.

\end{itemize}

% Generalized Turán and Generalized Rainbow Turán references

	\paragraph{Connections to generalized Turán and generalized rainbow Turán problems.}
	The method developed in this paper also relates to the broader framework of generalized Turán problems.
	The {\it generalized Turán number} $\ex(n, H, F)$ is the maximum number of copies of a fixed graph $H$ in an $n$-vertex $F$-free graph, first systematically studied by Alon and Shikhelman~\cite{AS16}. Recent work has obtained exact results for various specific graphs, including paths~\cite{HHL23} and their suspensions~\cite{HHM26}.
	In the context of rainbow subgraphs, the {\it generalized rainbow Turán number} $\exstar(n, H, F)$ is the maximum number of  copies of $H$ in a properly edge-colored $n$-vertex graph with no rainbow $F$.
	Gerbner, M\'esz\'aros, Methuku, and Palmer~\cite{GMMP21} initiated the study of generalized rainbow Turán numbers and established connections between the generalized Turán number and the generalized rainbow Turán number, analogous to the relationship between $\ex(n, H)$ and $\exstar(n, H)$.
	These generalizations offer a natural direction for future work: our heterochromatic-triangle-counting Lemma~\ref{lem:key} bounds the total number of edges in rainbow $H$-free graphs by constraining the number of triangles.
%	In particular, Hou, Hei, and Ma~\cite{HHM26} determined $\ex(n, K_3, H)$ where $H$ is the suspension of $P_5$, complementing the path results of~\cite{HHL23} and demonstrating the role of suspension graphs in generalized Tur\'an theory.
		A natural extension would be to bound the number of rainbow copies of a subgraph $H_0$ (not necessarily triangles) in a properly edge-colored graph that contains no rainbow $H$, thereby obtaining generalized rainbow Turán results for fans and wheels.

	% 参考文献列表（按原文编号，格式符合期刊规范）
\bibliographystyle{plainnat} % 数字标号引用风格（期刊常用）

\end{document}